\documentclass[11pt]{article}
\usepackage{amssymb}
\usepackage{amsmath}
\usepackage{latexsym}
\usepackage{amsthm}
\usepackage{mathptmx}

\makeatletter
\def\thmhead@plain#1#2#3{%
  \thmname{#1}\thmnumber{\@ifnotempty{#1}{ }#2}%
  \thmnote{ \the\thm@notefont(#3)}}
\let\thmhead\thmhead@plain
\def\swappedhead#1#2#3{%
  \thmnumber{#2}\thmname{\@ifnotempty{#2}{. }#1}%
  \thmnote{ \the\thm@notefont(#3)}}
\makeatother

\theoremstyle{definition} 

 \newtheorem{definition}{Definition}[section]
 \newtheorem{remark}[definition]{Remark}

\theoremstyle{plain}      

 \newtheorem{proposition}[definition]{Proposition}
 
 \newtheorem{theorem}[definition]{Theorem}
 
 \newtheorem{lemma}[definition]{Lemma}
 \newtheorem{defi}[definition]{Definition}
 
 \def \dem{\noindent{\sc Proof.~}}
\def \findem {\hfill{\hbox {\vrule\vbox{\hrule width 6pt\vskip 6pt\hrule}\vrule}}}

\newcommand{\cO}{{\mathcal O}}
\newcommand{\mm}{{\mathcal M}}
\def\OC{{\mathcal{O}}}
\def\O{{\mathcal{O}}}

\def\FC{{\mathcal{F}}}

\def\g{{\mathfrak{g}}}
\def\h{{\mathfrak{h}}}
\def\d{{\mathfrak{d}}}
\def\a{{\mathfrak{a}}}

\def\m{{\mathfrak{m}}}

\def\UM{{\mathbb{U}}}
\def\RM{{\mathbb{R}}}

\def\KM{{\mathbb{K}}}

\newcommand{\poly}{\operatorname{poly}}

\def \Vect {\mathop{\hbox{\rm Vect}}\nolimits}

\def \Alt {\mathop{\hbox{\rm Alt}}\nolimits}

\def \mod {\mathop{\hbox{\rm mod}}\nolimits}

\def \op {{\scriptstyle{\rm op}}}

\def \loc {{\scriptstyle{\rm loc}}}

\def \inv {{\scriptstyle{\rm inv}}}

\def \id {\mathop{\hbox{\rm id}}\nolimits}

\def \ve{\Lambda}

\newcommand{\pr}{\operatorname{pr}}

\def \Ob {\mathop{\hbox{\rm Ob}}\nolimits}
\def \Vect {\mathop{\hbox{\rm Vect}}\nolimits}

\def \Mer {\mathop{\hbox{\rm Mer}}\nolimits}

\def \Alt {\mathop{\hbox{\rm Alt}}\nolimits}
\def \CYB {\mathop{\hbox{\rm CYB}}\nolimits}
\def \DGLBA {\mathop{\hbox{\rm DGLBA}}\nolimits}

\def \DGQUE {\mathop{\hbox{\rm DGQUE}}\nolimits}

\def \Alt {\mathop{\hbox{\rm Alt}}\nolimits}

\def \Vect {\mathop{\hbox{\rm Vect}}\nolimits}

\def \d {\mathop{\hbox{\sl d}}\nolimits}
\def \cH {\mathop{\hbox{\rm cH}}\nolimits}
\def \mod {\mathop{\hbox{\rm mod}}\nolimits}

\newcommand{\Lie}{\operatorname{Lie}} 

\newcommand{\alg}{\operatorname{alg}} 
\newcommand{\G}{\operatorname{G}} 
\newcommand{\Gi}{\G_\infty} 
\newcommand{\Ger}{\operatorname{Ger}}

\def \ve {\Lambda}
\def \g {\mathfrak{g}}
\def \a {\mathfrak{a}}

\def \int {{\scriptstyle{\rm int}}}

\newcommand{\cT}{{{}^cT}}
\newcommand{\uTTU}{{\underline{{}^cT}T_+U}}
\newcommand{\uTTH}{{\underline{{}^cT}T_+H}}
\newcommand{\CuTTU}{{C(\underline{{}^cT}T_+U)}}
\newcommand{\CuTTUt}{{C(\underline{{}^cT}T_+\tilde{U})}}
\newcommand{\uTCuTTU}{{\underline{{}^cT}C(\underline{{}^cT}T_+U)}}

\newcommand{\uTA}{{\underline{{}^cT}A}}
\newcommand{\uTCuTA}{{\underline{{}^cT}C(\underline{{}^cT}A)}}

\newcommand{\uTE}{{\underline{{}^cT}E}}
\newcommand{\uTV}{{\underline{{}^cT}V}}
\newcommand{\uTpE}{{\underline{{}^cT}(E)}}

\newcommand{\dcH}{b_{\cH}}

\newcommand{\Tpi}{T_{\poly}^{\inv}}

\newcommand{\Dpi}{D_{\poly}^{\inv}}
\newcommand{\fa}{\varphi_{\alg}}
\newcommand{\fg}{\varphi_{{\Ger}_\infty}}
\newcommand{\fl}{\varphi_{\Lie}}
\newcommand{\fld}{\varphi_{\Lie}'}
\newcommand{\fgi}{\varphi_{\Gi}}

\newcommand{\DQ}{\operatorname{DQ}}
\newcommand{\LG}{\operatorname{L-G}}
\newcommand{\LGi}{\LG_\infty}
\newcommand{\Uhg}{{U_\hbar({\a})}}

\begin{document}

\title{Quantization of $r-Z$-quasi-Poisson manifolds and related modified classical dynamical $r$-matrices}

\author{Gilles Halbout\\
{\small{}}\cr
{\small{ \! Institut de Math\'ematiques et de Mod\'elisation de Montpellier}}\cr
{\small{UMR 5149 de l'Universit\'e Montpellier II et du CNRS}}\cr
{\small{ CC 5149,  Place Eug\`ene Bataillon }}\cr
{\small{F-34 095 Montpellier CEDEX 5}}\cr
{\small{ \! e-mail:\,\texttt{ghalbout@math.univ-montp2.fr}}}\cr
}
\markboth{Gilles Halbout}
{Quantization of $r-Z$-quasi-Poisson manifolds and related modified classical dynamical $r$-matrices}


\maketitle

\abstract{{\small 
Le $X$ be a $C^\infty$-manifold and $\g$ be a finite dimensional
Lie algebra acting freely on $X$. Let $r \in \ve^2(\g)$ be such that
$Z=[r,r] \in \ve^3(\g)^\g$. In this paper we prove that 
every quasi-Poisson $(\g,Z)$-manifold can be quantized.
This is a generalization of the existence of a twist
quantization of coboundary Lie bialgebras (\cite{EH})
in the case $X=G$ (where $G$ is the simply connected
Lie group corresponding to $\g$). 
We deduce our result from a generalized formality theorem.
In the case $Z=0$, we get a new proof of the existence of (equivariant) formality theorem and so (equivariant) quantization of Poisson manifold
({\it cf.} \cite{Ko,Do}). As a consequence of our results, we get quantization of modified classical dynamical $r$-matrices
over abelian bases in the reductive case.}}
\vskip20pt

\centerline {\bf $0$. Introduction }

\vskip20pt

Throughout this paper, the ground field will be $\RM$.
Let $\g$ be a finite dimensional Lie algebra with a fixed element
$r \in \ve^3(\g)$ such that 
$[r,r]=Z \in \ve^3(\g)^\g$. In \cite{AK,AKM}, quasi-Poisson manifolds 
were introduced as a generalization of Poisson $\g$-manifolds with
Poisson bracket satisfying the Jacobi identity up to an invariant trivector
corresponding to $Z$. More precisely~:
\begin{defi}
A quasi-Poisson $(\g,Z)$-manifold is a $\g$-manifold $X$ with an invariant bivector
$\pi$ such that the Schouten bracket $[\pi,\pi]_S$ equals
$\gamma^{\otimes 3}(Z)$, where $\gamma$~:
$\g \to \Vect(X)$ is the action homomorphism.
\end{defi}
\noindent 
The Schouten bracket will be descibed later. Thus the Poisson bracket $\{-,-\}$ associated to $\pi$ satisfies
$$\{\{f,g\},h\}+\{\{g,h\},f\}+\{\{h,f\},g\}=m_0(\gamma^{\otimes 3}(Z)(f \otimes g
\otimes h)),$$
where $m_0$ is the usual multiplication.
In the framework of deformation quantization (see
\cite{BFFLS1, BFFLS2}), Enriquez and Etingof defined 
the quantization of quasi-Poisson
ma\-nifolds 
in \cite{EE1}~: let $\hbar$ be a formal parameter and
$\Phi=1+\frac{\hbar^2}{6} Z + O(\hbar^3) 
\in (U(\g)^{\otimes 3})^\g
[[\hbar]]$ be an associator for $\g$ (Drinfeld proved in \cite{Dr}, Proposition 3.10, that such an associator always exists). 
\begin{defi}
\label{defstar}
A quantization of $X$ associated to $\Phi$ is an invariant star-product 
$\star$ on $X$,i.e. an invariant bidifferential operator
on $C^\infty(X)$,
which
satisfies $f\star g= fg +O(\hbar)$ and the equation
$$f \star g - g \star f = \hbar \{ f,g\} + O (\hbar^2),$$
and is associative in the tensor category of
$(U(\g)[[\hbar]],\Phi)$-modules. This means,
$$m_\star(m_\star \otimes 1)=m_\star (1 \otimes m_\star) \gamma^{\otimes 3}(\Phi),$$
on $C^{\infty}(X)^{\otimes 3}$, where $m_\star(f
\otimes g)=f \star g$.
\end{defi}
\noindent They also conjectured that such quantizations always exist when the action of $\g$ 
on the
quasi-Poisson manifold
$X$ is free. Note that when the action is not free, 
Fronsdal (\cite{Fr}) gave in 1978 a counter-example where
such quantization is impossible even in the symplectic case.
From now on, we will suppose that the manifold $X$ is a $G$-bundle over
a manifold $M$, where $G$ is the simply connected Lie group corresponding to $\g$.
In the case $G=\{e\}$, the conjecture is equivalent to the existence
of star-products and was proved by Kontsevich (\cite{Ko}).
In the case $Z=0$, the conjecture
follows from the equivariant formality theorem of Dolgushev (\cite{Do}).

\medskip

In the general case,
$\gamma^{\otimes 3}(Z)$ commutes with all the left invariant
polyvector fields in the following sense~:
\begin{equation}
\label{cond1}
[\gamma^{\otimes 3}(Z),X]_S=0,~\hbox{ for all invariant 
polyvector fields } X.
\end{equation}
Moreover, for $\Phi$ an associator, $\gamma^{\otimes 3}(\Phi)$ commutes
with all the invariant differential operators
in the following sense:
\begin{equation}
\label{cond2}
[\gamma^{\otimes 3}(\Phi),C]_G= 0,\hbox{ for all invariant differential
operator }C
\end{equation}
(the Gerstenhaber bracket $[-,-]_G$ will be described later in this paper).
From now on, if $s \in \ve^k(\g)$, we will denote $s$ instead of $\gamma^{\otimes k}(s)$ when no
confusion is possible.

In this paper, we prove that there exists (a least) one associator for $\g$ such that Enriquez-Etingof's conjecture is true:
\begin{theorem}
\label{theoprinc}
Let $r \in \ve^3(\g)$ such that 
$[r,r]=Z \in \ve^3(\g)^\g$. There exists
$\Phi=1+\frac{\hbar^2}{6} Z + O(\hbar^3) 
\in (U(\g)^{\otimes 3})^\g
[[\hbar]]$ and a deformation $\g_\hbar$ of the
Lie algebra $\g$ such that for every invariant
bivector $\pi$ satisfying $[\pi,\pi]_S=\gamma^{\otimes 3}(Z)$,
the quasi-Poisson manifold $(X,\pi)$ admits a quantization associated to
$(\Phi,\g_\hbar)$ i.e. a multiplication associative in the tensor category of
$(U(\g_\hbar)[[\hbar]],\Phi)$-modules.
\end{theorem}

\smallskip

To prove this theorem, we will construct a formality 
between invariant polyvector and polydifferential operator as stated in Theorem
\ref{theo:linfprinc}. We first prove a local version of this theorem in the 
case $X=\RM^n \times \g$.
Using Fedosov's resolutions we will be able to get a
global version. We then get the wanted invariant
star-product on the manifold $X$ and classification of such deformations. We will then discuss the relation with
quantization of modified classical dynamical $r$-matrices.

\begin{remark}
As a particular case, our results give a new proof of Kontsevich (and Dolgushev
for equivariant) formality theorem. One can see this approach as related to Merkulov's work (see \cite{Me}) for quantization of Lie bialgebras. In our
work the use of a graded version of Etingof-Kazhdan theorem was a crucial step
to go from quantization of Lie bialgebra to quantization of Poisson manifolds.
\end{remark}

\bigskip

The paper is organized as follows:

\noindent - In Section \ref{sec:structure}, we recall definitions of $L_\infty$-structures and formality morphisms. 

\noindent - In Section \ref{sec:EK}, we give a graged version of  quantization of Lie bialgebras: in
particular, we get differential graded Etingof-Kazhdan quantization/dequantization functors.

\noindent - In Sections \ref{sec:funct} and \ref{sec:reso}, we construct two useful functors between Lie and Gerstenhaber algebras ``up to homotopy'' and prove the existence of two resolutions for those algebras.

\noindent - In Section \ref{sec:Ha}, we prove the existence of $L_\infty$-morphisms between DG Lie bialgebras and the Gerstenhaber algebra of their Etingof-Kazhdan quantization  

\noindent - In Section \ref{sec:phi}, we transpose the algebra structures into the category of 
$(U(\g)[[\hbar]],\Phi)$-modules. We define the graded Lie bialgebra $\tilde{\g}=\RM\oplus
V[1] \oplus V^* \oplus \g$, the direct sum of the Eisenberg Lie
algebra  $E=\RM\oplus V[1] \oplus V^*$ and the Lie bialgebra $(\g,[r,-])$ which corresponds locally to the algebra of invariant poly-vectors. We 
prove the existence of the local wanted $L_\infty$-morphism.

\noindent - In Section \ref{sec:glob}, we show that this $L_\infty$-morphism can be globalized and prove our main theorem.

\noindent - In Section \ref{sec:discus}, we discuss relation between our quantization and quantization of modified classical dynamical $r$-matrices.

\medskip

\subsection*{Notations}

We use
the standard notation for the coproduct-insertion maps:
we say that an ordered set is a pair of a finite set $S$ and a 
bijection $\{1,\ldots,|S|\} \to S$. 
For $I_1,\dots,I_m$ disjoint ordered subsets of $\{1,\dots,n\}$,
$(U,\Delta)$ a Hopf algebra and $a \in  U^{\otimes m}$,
we define
$$a^{I_1,\dots,I_m}= \sigma_{I_1,\ldots,I_m} \circ 
(\Delta^{(|I_1|)}\otimes \cdots \otimes \Delta^{(|I_m|)})(a),
$$
with $\Delta^{(1)}=\id$, $\Delta^{(2)}=\Delta$, 
$\Delta^{(n+1)}=({\id}^{\otimes n-1} \otimes \Delta)\circ \Delta^{(n)}$, 
and $\sigma_{I_1,\ldots,I_m} : U^{\otimes \sum_i |I_i|} \to 
U^{\otimes n}$ is the 
morphism corresponding to the map $\{1,\ldots,\sum_i |I_i|\} 
\to \{1,\ldots,n\}$ taking $(1,\ldots,|I_1|)$ to $I_1$, 
$(|I_1| + 1,\ldots,|I_1| + |I_2|)$ to $I_2$, etc. 
When $U$ is cocommutative, this definition depends only on
the sets underlying $I_1,\ldots,I_m$.  

\smallskip

Until the end of this paper, although we will often omit to mention it, we will always deal with graded structures.

\medskip

\subsection*{Acknowledgements}
I would like to thank D. Calaque and B. Enriquez for discussions. Part of this work
was done during the program ``Poisson Sigma Models, Lie Algebroids, Deformations, and Higher Analogues'' at the Erwin Schr\"odinger International
Institute for Mathematical Physics.

\medskip

\section{$L_\infty$-structures}\label{sec:structure}

\medskip

\subsection{Definitions}

\medskip

Let us recall definitions of $L_\infty$-algebras and $L_\infty$-morphisms.
Let $A$ be a graded vector space.
We denote $T_+A=T_+(A[-1])$ the free tensor algebra (without unit) of $A$
which, equipped with the coshuffle coproduct, is a bialgebra.
We also denote $C(A)=S(A[-1])$ the free graded commutative algebra
generated by $A[-1]$, seen as a quotient
of $T_+A$. The coshuffle coproduct is still well defined on $C(A)$ which 
becomes a cofree cocommutative coalgebra on $A[-1]$.
We also denote $\Lambda A=S(A[1])$, the analogous graded commutative algebra
generated by $A[1]$ (in particular, for $A_1,A_2 \in A$, $A_1 \Lambda A_2$ stands for the
corresponding quotient of $A_1[1] \otimes A_2[1]$ in $\Lambda A$).
We will use the notations $T_+^nA$, $\Lambda^nA$ and $C^n(A)$ for the elements
of degree $n$.

\begin{definition}\label{DefinitionL_infini}
A vector space $A$ is endowed with a $L_\infty$-algebra 
(Lie algebra ``up to homotopy'') structure
if there are degree one linear  maps $d^{1,\dots,1}$: 
$\Lambda^kA \rightarrow A[1]$ such that 
the associated coderivations 
(extended with respect to the
cofree cocommutative structure on $\Lambda A$) d: $\Lambda A \rightarrow \Lambda A$, 
satisfy $d \circ d=0$ where $d$ is the coderivation 
$$d=d^1+d^{1,1}+\cdots+d^{1,\dots,1}+\cdots.$$
\end{definition}
\noindent In particular, a differential Lie algebra $(A,b,[-,-])$ is
a $L_\infty$-algebra with structure maps $d^1=b[1]$, $d^{1,1}=[-,-][1]$ and
$d^{1,\dots,1}$: 
$\Lambda^k A \rightarrow A[1]$ are $0$ for $k \geq 3$.
One can now define the generalization of Lie algebra morphisms:
\begin{definition}\label{Linfini}
A $L_{\infty}$-morphism between two $L_\infty$-algebras 
$(A_1,d_1=d^1_1+\cdots)$ and $(A_2,d_2=d^1_2+\cdots)$ is 
a morphism of codifferential cofree coalgebras, of degree $0$,
$$
\varphi: ~(\Lambda A_1,d_1) \rightarrow (\Lambda A_2,d_2).
$$
\end{definition}
\noindent In particular 
$\varphi \circ d_1 = d_2 \circ \varphi $. 
As $\varphi$ is  a morphism of cofree cocommutative  coalgebras, $\varphi$
is determined by its image on the cogenerators, i.e., by its components:
$\varphi^{1,\dots,1}: \Lambda^k A_1\to A_2[1].$ 

\medskip

Let $E$ be a graded vector space.
Let us denote $\cT (E)$ the cofree tensor coalgebra of $E$
with coproduct $\Delta'$. Equipped with the
shuffle product $\bullet$ (defined on the cogenerators  $\cT (E) \otimes  \cT (E) \to E$
as $\pr \otimes \varepsilon + \varepsilon \otimes \pr$, where $\pr$ :  $\cT (E) \to E$ is
the projection and $\varepsilon$ is the counit), it is a bialgebra. 
Let $\cT_+ (E)$ be the augmentation ideal.
We denote $\uTpE= \cT_+ (E)/(\cT_+ (E) \bullet \cT_+ (E))$ the quotient by the shuffles.
It has a graded cofree Lie coalgebra structure (with coproduct $\delta = \Delta' -{\Delta'}^\op$).
Then $S(\uTpE[1])$ has a structure of cofree coGerstenhaber algebra (i.e. equipped with
cofree coLie and cofree cocommutative coproducts satisfying compatibility condition).
We use the notation $\underline{\cT^n} (E)$ for the elements of degree $n$.

\begin{remark}
One could also define $G_\infty$-structures. Most of the $L_\infty$-morphism constructed in this paper are also $G_\infty$-morphisms between corresponding
$G_\infty$-structures. Definitions and extensions to $G_\infty$-structures
can be found in \cite{Ha}.
\end{remark}

\smallskip

\section{Etingof-Kazhdan functors}\label{sec:EK}

\medskip

\subsection{QUE and QFSH algebras} \label{subsec:duality}

\medskip

We recall some facts from \cite{Dr} (proofs and definitions can be found in
\cite{Gav}). Let us denote by ${\bf QUE}$ the category of quantized universal
enveloping (QUE) algebras and by ${\bf QFSH}$ the  category of quantized formal
series Hopf (QFSH) algebras. 
Let us recall the definition of FSH and QFSH algebras:
\begin{definition}
A FSH algebra is a Hopf algebra of power series  isomorphic as an algebra to 
$\KM[[\{u_i | i \in J\} ]]$ (for some set $J$).
\end{definition}
There is an equivalence of categories between the category of FSH algebra and the 
category of Lie coalgebra (LC algebra), sending $\O_\h$ to $\h={\O_\h}_+ \slash {\O_\h}_+^2$
where ${\O_\h}_+$ is the maximal ideal of ${\O_\h}$.
\begin{definition}
A QFSH algebra is a Hopf algebra $H$, which is a topologically free  $\KM[[\hbar]]$-module, 
such that $H_0:=H/\hbar H$ is isomorphic to a FSH algebra
\end{definition}

\smallskip

Let us give an example of a FSH algebra, very important in this paper: 
let $V$ be a vector space and ${\cT}V $ (defined in the previous section)
the cofree coalgebra, equipped with
the shuffle product. 
Let us now complete $\cT V$. The algebra $\cT V$ is a graded algebra with
$V$ being the set of elements of degree $1$. Let us denote $\mm_{\cT V}$ the set
of elements of degree $\geq 1$. Finally, we denote $\widehat{\cT} V$ the
commutative cofree bialgebra, $\mm_{\cT V}$-adic
completion of $\cT V$.
\begin{proposition}\cite{Ha}
\label{FSHA}
$\widehat{\cT} V$ is the FSH algebra $\O_{\uTV}$ associated with the
Lie coalgebra $\uTV=\cT_+ V/(\cT_+ V)^2$, which
is the cofree Lie coalgebra over $V$.
\end{proposition}



We have covariant functors ${\bf QUE} \to {\bf QFSH}$, $U\mapsto U'$
and ${\bf QFSH} \to {\bf QUE}$, $\cO\mapsto \cO^\vee$. These functors are
also inverse to each other. 

$U'$ is a subalgebra of $U$
defined as follows:
for any ordered subset 
$ \, \Sigma = \{i_1, \dots, i_k\} \subseteq \{1, \dots, n\} \, $ 
with  $ \, i_1 < \dots < i_k \, $,  \, define the morphism 
$ \; j_{\scriptscriptstyle \Sigma} : U^{\otimes k} \longrightarrow
U^{\otimes n} \; $  by  $ \; j_{\scriptscriptstyle \Sigma}
(a_1 \otimes \cdots \otimes a_k) := b_1 \otimes \cdots \otimes
b_n \; $  with  $ \, b_i := 1 \, $  if  $ \, i \notin \Sigma \, $  and 
$ \, b_{i_m} := a_m \, $  for  $ \, 1 \leq m \leq k $;  then set 
$ \; \Delta_\Sigma := j_{\scriptscriptstyle \Sigma} \circ \Delta^{(k)}
\, $,  $ \, \Delta_\emptyset := \Delta^{(0)} \, $,  and
 $$  \delta_\Sigma := \sum_{\Sigma' \subset \Sigma} {(-1)}^{n- \left|
\Sigma' \right|} \Delta_{\Sigma'}  \; ,  \qquad \quad  \delta_\emptyset
:= \epsilon \; .   $$  
We shall also use the notation  $ \, \delta^{(n)} := \delta_{\{1, 2,
\dots, n\}} \, $,  $ \, \delta^{(0)}
:= \delta_\emptyset \, $,  and the useful formula  
  $$  \delta^{(n)} = {({id}_{\scriptscriptstyle U} - \epsilon)}^{\otimes n}
\circ \Delta^{(n)} \, .  $$  
Finally, we define
  $$  U' := \big\{\, a \in U \,\big\vert\, \delta^{(n)}(a) \in
h^n U^{\otimes n} \, \big\} \quad ( \subseteq U )  $$
and endow it with the induced topology.  

\medskip

On the other way, $\cO^\vee$
is the $\hbar$-adic completion of $\sum_{k\geq 0} \hbar^{-k} \mm^k \subset 
\cO[1/\hbar]$ (here $\mm\subset \cO$ is the maximal ideal).

\medskip

\subsection{The functor $\DQ$}

\medskip

In \cite{GH}, a generalization of Etingof-Kazhdan theorem (\cite{EK}) was proved
in an appendix by Enriquez and Etingof:
\begin{theorem}
\label{thm:dqdg}
We have an equivalence of categories
$$DQ_\Phi~:~\DGQUE \to \DGLBA_h$$
from the category of differential graded  
quantized universal enveloping super-algebras to that of differential graded Lie
super-bialgebras such that if
$U \in \Ob(\DGQUE)$ and $\a=DQ(U)$, then
$U/hU=\UM(\a/h\a)$, where $\UM$ is the universal algebra functor, taking a
differential graded Lie super-algebra to a differential graded super-Hopf 
algebra.
\end{theorem}
Here $\Phi$ is a Drinfeld associator. We will use any of these functors and denote it $\DQ$.

\medskip

%

\section{Two functors}\label{sec:funct}

\medskip

\subsection{Functor $\LG$} Let $(\h,\delta,d)$ be a differential Lie bialgebra. Let $C(\h)=S(\h[-1])$ be the free graded
commutative algebra generated by $\h$. Recall from the previous subsection that $C(\h)$ is also a cofree coalgebra and that 
coderivations $C(\h) \to C(\h)$ are defined by their images in $\h$. Thus, one
easily checks that the coderivation $[-,-]$: $C(\h) \to C(\h)$ extending the Lie
bracket (with degree shifted by one) defines a Lie (even Gerstenhaber) algebra structure
on $C(\h)$. Moreover,
one can extend maps $d$: $\h \to \h$ and $\delta$: $\h \to S^2(\h[-1])$ on
the free commutative algebra $C(\h)$ so that
$(C(\h),[-,-],\wedge,d+\delta)$ is a differential Gerstenhaber algebra. The
differential $\delta$ is actually the Chevalley Eilenberg differential: the space
$C(\h)=S^*(\h[-1])$ is isomorphic to the standard complex $(\Lambda^*(\h))[-*]$ and
$\delta$ is simply the differential given by 
the underlying Lie coalgebra structure of $\h$.
\begin{proposition}\cite{Ha}
Any DGLA morphism $f:\h_1\to \h_2$ can be 
extended into a DGLA (and even differential graded Gerstenhaber) morphism $C(f):C(\h_1) \to C(\h_2)$ 
of free commutative algebras. 
This defines an exact functor ${\LG}$ 
from differential Lie bialgebras to differential Gerstenhaber algebras which sends 
$\h$ to $C(\h)$. Quasi-isomorphisms $(\h_1,d_1)\to (\h_2,d_2)$ 
induce a quasi-isomorphisms $(C(\h_1),d_1,\delta_1)\to (C(\h_2),d_2,\delta_2)$.
\end{proposition}

\medskip

\subsection{Functor $\LGi$} Consider now the category ${\rm CFDLB}$ of  
differential Lie bialgebras which are cofree as a Lie coalgebra. 
In other words we are interested in cofree Lie coalgebra $\uTpE$ 
on a graded vector space $E$ together with a differential $\ell$ 
and a cobracket $L$ on $\uTpE $ that makes it a differential Lie bialgebra. 
As $\uTpE$ is cofree, the differential is uniquely determined 
by its restriction to cogenerators $l^p$: $\underline{\cT^p}(E) \to E$. 
Similarly, the Lie bracket is uniquely determined by maps $L^{p_1,p_2}$:  
$\underline{\cT^{p_1}}(E)\Lambda \underline{\cT^{p_2}}(E) \to E$.
\begin{proposition}\cite{Ha}
Restriction map $\underline{\cT^p}(E) \to E$ defines an exact functor $\LGi$
from ${\rm CFDLB}$ to the category of $G_\infty$ (and so Lie)-algebras.
\end{proposition} 

\medskip

Until the end of the paper, we will use the notations $TE$ for $T(E[-1])$ 
and $\uTE$ 
for $\underline{T}(E[1])$.

\smallskip

\section{Two resolutions}\label{sec:reso}

\medskip

\subsection{bialgebra structure on $\cT T_+U$}

Here, we will define a bialgebra structure on $\cT T_+U$.
One can construct a bialgebra structure on the space of
Hochschild co\-chains of an algebra using the brace operations.
In our case, we will firstly generalize the definition of brace operations for a general Hopf algebra.
More precisely, let $(H,\Delta_\hbar,\times)$ be a Hopf algebra (in our case $H$ will be
the Etingof-Kazhdan quantization $\Uhg$ of the Lie bialgebra $\a$).
We will define a brace structure on the cofree tensor coalgebra $\cT T_+H$ of
the free tensor algebra $T(H[-1])$ without unit. 
To distinguish the two tensor products, we denote $\otimes$ the
tensor product on $T_+H$ and $\boxtimes$ the tensor product
on $\cT T_+H$.
\begin{definition}
\label{brace}
We define brace operations on $\cT T_+H$ by extending the following maps
given on the cogenerators of the cofree coalgebra $\cT T_+H$:
\begin{enumerate}
\item $B^0=0$,
\item $B^1=\dcH$ (the coHochschild coboundary on $T_+H$),
\item $B^2$ : $\alpha \boxtimes \beta \mapsto \alpha \otimes \beta$,
\item $B^n=0$  for $n >2$,
\item $B^{0,1}=B^{1,0}=\id$,
\item $B^{0,n}=B^{n,0}=0$ for $n \geq 1$,
\item $B^{1,n}$ : $(\alpha,\beta_1\boxtimes \cdots \boxtimes \beta_n) \mapsto$ 
\begin{multline*}
\sum_{\stackrel{0 \leq i_1,\dots,i_m+k_m\leq n}
{i_l+k_l\leq i_{l+1}}} (-1)^\varepsilon \alpha^{1,\dots,i_1+1\cdots i_1+k_1,\dots,i_m+1 \cdots i_m+k_m,
\dots,n} \times \\
1^{\otimes i_1} \otimes \beta_1 \otimes 1^{\otimes i_2-(i_1+k_1)}\otimes \beta_2 \otimes \dots \otimes \beta_n \otimes
1^{\otimes n-(i_m+k_m)},
\end{multline*}
where $k_s=|\beta_s|$, $n=|\alpha|+\sum_sk_s - m$ and $\varepsilon = \sum_s (k_s-1)i_s$,
\item $B^{k,l}=0$ for $k >1$.
\end{enumerate}
Operations $(2),(3)$ and $(4)$ define a differential $d$ and $(5),(6),(7)$ and $(8)$ define a product
$\star$ deforming the shuffle product.
\end{definition}
Note that, when $H=U(\a)$, the enveloping
algebra of a Lie algebra $\a$, $T(H[-1])$ can be seen as the space of
invariant polydifferential operators over the Lie group corresponding to $\a$
and in that case, our definition coincides with usual braces operations.

We have:
\begin{theorem}\cite{Ha}
\label{bial}
The brace operations of Definition \ref{brace} define a differential bialgebra
structure on the cofree tensor coalgebra $\cT T_+H$, with product  $\star$
extending $\sum B^{p_1,p_2}$ and differential $d$ extending
$\sum B^p$.
\end{theorem}

\medskip

Let us now complete $\cT T_+H$ as in section \ref{sec:EK} with 
$V=T_+H$. We get  a commutative cofree bialgebra $\widehat{\cT} T_+H$, the  $\mm_{\cT T_+H}$-adic
completion of $\cT T_+H$ (where $\mm_{\cT T_+H}$ is the
maximal ideal of $\cT T_+H$).
Let us consider the free $\KM[[\nu]]$-module $\widehat{\cT} T_+H[[\nu]]$.
One can now replace the operations $B^{p,q}$ of Definition \ref{brace}
with $\KM[[\nu]]$-linear operations $\nu^{p+q-1}B^{p,q}$. Those operations 
are well defined on the completion $\widehat{\cT} T_+H[[\nu]]$ as this space is complete
for the grading induced by the degree in $\cT T_+H=\cT V$
plus the $\hbar$-adic valuation and because
the operations we just defined are homogeneous for this grading.
Thus we get a morphism of differential bialgebra
\begin{align}
\label{formule:iso}
I_\nu:~({\cT} T_+H, \star,\Delta,d) & \to ({\cT} T_+H[[\nu]][\nu^{-1}],\star_\nu,\Delta_\nu,d_\nu)   \cr
x & \mapsto \nu^{-|x|} x, 
\end{align}
where $|x|$ is the degree in ${\cT}$.
The morphism $I_\nu$ extends to 
$I_\nu:$ $({\cT} T_+H[[\nu]], \star,\Delta,d)  \to ({\cT} T_+H[[\nu]][\nu^{-1}],\star_\nu,\Delta_\nu,d_\nu)$
which restricts to
\begin{equation}
\label{formuleprime}
I_\nu':~(\widehat{\oplus}\nu^n{\cT}^n T_+H[[\nu]], \star,\Delta,d)  \to (\widehat{\cT} T_+H[[\nu]][\nu^{-1}],\star_\nu,\Delta_\nu,d_\nu)  
\end{equation}

We have:
\begin{proposition}\cite{Ha}
\label{pro:QF}
The algebra $(\widehat{\cT} T_+H[[\nu]],\star_\nu,\Delta_\nu,d_\nu)$ is a QFSHA.
The underlying differential Lie bialgebra structure
on $\uTTH$
 is given by the Gerstenhaber bracket 
$$[\alpha,\beta]_G=B^{1,1}(\alpha,\beta) - (-1)^{(|\alpha|-1)(|\beta|-1)}B^{1,1}(\beta,\alpha) $$
and coHochschild differential
$$\dcH(\alpha)=[1\otimes 1, \alpha]_G,$$
for $\alpha,\beta \in TH$ and then naturally extended on $\uTTH$ using the cofree Lie cobracket.
\end{proposition}

\smallskip

\begin{remark}
Let now $H$ be the QUE algebra $U=U_\hbar(\a)$. 
We have proved that 
$T_+U$ can be equipped with a $G_\infty$-structure.
Since the cofree Lie coalgebras are rigid, the differential Lie bialgebra corresponding to $\widehat{\cT} T_+U[[\nu]]$ through
Etingof-Kazhdan dequantization functor $\DQ$ is isomorphic to
$\uTTU[[\nu]]$ as a  $\KM[[\nu]]$-Lie coalgebra, and is therefore free.
\end{remark}

\medskip

\subsection{A bialgebra quasi-isomorphism $\fa$ : $U\to (\widehat{\cT} T_+U)^\vee$}

\medskip

We have:
\begin{proposition}
Let $U$ be a QUE algebra.
One can define a bialgebra quasi-isomorphism
$\phi_{\alg}$: $U\to \widehat{\cT} T_+U$ from the
bialgebra $(U,\Delta_\hbar,\times)$ to
the bialgebra $({\cT} T_+U,\Delta,\star)$ whose
structure was described in the previous
section. 
\end{proposition}

\medskip

Let $U' \subset U$ (see section \ref{sec:EK}). 
\begin{proposition}\cite{Ha}
We have a bialgebra quasi-isomorphism
$\fa$ : $(U',\times)\to (\widehat{\cT} T_+U,\star_\hbar)$ of QFSH algebra,
where  $(\widehat{\cT} T_+U,\star_\hbar)$ is
$(\widehat{\cT} T_+U[[\nu]],$ $\star_\nu)\slash (\nu=\hbar)$
($\widehat{\cT} T_+U[[\nu]]$ is the free $\KM[[\hbar]]$-module
 defined in the previous section: we 
 d the operations $B^{p,q}$ into $\nu^{p+q-1}B^{p,q}$).
\end{proposition}
 
\smallskip

Finally, applying to $\fa$ the derived Drinfeld functor $(-)^\vee$, we get 
a bialgebra quasi-isomorphism
$\fa$: $U\to (\widehat{\cT} T_+U)^\vee$.

\medskip

\subsection{A Lie bialgebra quasi-isomorphism $\fl$ : $\uTA \to \uTCuTA$}

\medskip

Let $A$ be a vector space.
Suppose now that the cofree Lie coalgebra $\uTA$
has a structure $(\uTA,\delta,[-,-],d)$ of a differential Lie bialgebra.
Using the functor ${\LG}$ (see section \ref{sec:funct}), one gets a differential Gerstenhaber algebra
$(C(\uTA),[-,-],\wedge,d+\delta)$.
One can extend the structure maps on
the cofree Lie coalgebra $\uTCuTA$
and one gets 
a differential cofree Lie bialgebra 
$(\uTCuTA,\delta',[-,-],d+\delta+\wedge)$ (we will set $d^1=d + \delta$ and  $d^2=\wedge$).
\begin{proposition}\cite{Ha}
Let $\fl$ be the composition map $\fl= {\cT}i \circ \bar{\delta}$
of a map
\begin{equation*}
\begin{array}{rl}
\bar{\delta}: \uTA &\to \underline{\cT}(\uTA),\\
x & \mapsto x + \sum_{k \ge   2}\bar{\delta}_k(x),
\end{array}
\end{equation*}
where $\bar{\delta}_k$ is built using iterates of $\delta$, with
${\cT}i$:
$ \underline{\cT}(\uTA) \to \uTCuTA$ which is $ \underline{\cT}$ of the inclusion $i$: $\uTA[-1] \to C(\uTA)$.
Then $\fl$ is a differential Lie bialgebra quasi-isomorphism $\fl$ : $\uTA \to \uTCuTA$.
\end{proposition}

\medskip

\section{$L_\infty$-morphism for Lie bialgebras}\label{sec:Ha}

\medskip

\subsection{A Lie bialgebra quasi-isomorphism $\fld$: $\a\to\uTTU$}

\medskip

Let $(\a,\delta_\hbar)$ be a graded Lie bialgebra. We write
$\delta_\hbar=\hbar \delta_1 + \hbar^2 \delta_2 + \cdots$.
Let $(\Uhg,\Delta_\hbar)$ be the Etingof-Kazhdan canonical quantization
of $(\a,\delta_\hbar)$. 
We denote $U=\Uhg$ for short.
In section \ref{sec:reso}, we proved the existence of a bialgebra structure on $\cT T_+U$ and a
bialgebra quasi-isomorphism $\fa$: $U\to (\cT T_+U)^\vee$.
Thanks to Etingof-Kazhdan dequantization functor (see section
\ref{sec:EK}), and the fact
that $(\cT T_+U)^\vee$ is a QUE algebra quantizing $\uTTU$
(see section \ref{sec:reso}),
we get a Lie bialgebra quasi-isomorphism $\fld$: $\a\to\uTTU$.

\medskip

\subsection{Inversion of formality morphisms}\label{sec:inverse}

\medskip

Let us recall Theorem 4.4 of Kontsevich (\cite{Ko}):
\begin{theorem}\label{theo:koninv}
Let $\g_1$ and $\g_2$ be two $L_\infty$-algebras and $\FC$
be a $L_\infty$-morphism from $\g_1$ to $\g_2$.
Assume that
$\FC$ is a quasi-isomorphism. Then there exists an $L_\infty$-morphism
from $\g_2$ to $\g_1$ inducing the inverse isomorphism
between associated cohomology of complexes.
\end{theorem}

\begin{remark}
\label{greg}
We know the existence of a similar
$G_\infty$-version of this theorem. This result would imply the
existence of corresponding $G_\infty$-morphisms.
\end{remark}

\medskip

\subsection{$L_\infty$-morphism for Lie bialgebras}
 
\medskip

Let us summarize functors and quasi-isomorphisms constructed in
the previous sections in the following diagram:

\begin{equation*}
\begin{array}{ccccccccc}
\uTCuTTU & &\CuTTU[1] & = & \CuTTU[1]  &&\uTTU &&(\widehat{\cT} T_+U)^\vee\\[0.3cm]
\uparrow^{\fl} & \hskip-0.2cm\stackrel{{\LGi}}{\longrightarrow}\hskip-0.2cm& \uparrow^{\fgi} & & \uparrow^{\fg}&
\hskip-0.2cm\stackrel{{\LG}}{\longleftarrow}\hskip-0.2cm& \uparrow^{\fld}&\hskip-0.2cm\stackrel{{\DQ}}{\longleftarrow}& \hskip-0.2cm\uparrow^{\fa}\\[0.3cm]
\uTTU & &T_+U[1]  &  & C(\a)[1]  &&\a &&U=\Uhg.
\end{array}
\end{equation*}
Thus, thanks to section \ref{sec:inverse}, the composition $\varphi$: $C(\a) \to T_+U$ of $\fg$ with the inverse of $\fgi$ gives the wanted quasi-isomorphism.

\begin{theorem}\cite{Ha}
the map $\varphi$ $C(\a) \to T_+U$ is a $L_\infty$-quasi-isomorphism that
maps $v \in C(\a)$ to $\Alt(v) \in T_+U \mod \hbar$.
\end{theorem}

\medskip

\subsection{$L_\infty$-morphism for $X=\RM^n \times \g$}
\label{S5.4}
\medskip

We will now consider $X=\RM^n \times \g$ and $r\in \g \wedge \g$ such that $[r,r]=Z$. So $(\g,[r,-])$ is a Lie bialgebra. Let us set $V=\RM$.
From now on we will consider the graded Lie bialgebra $\tilde{\g}=\RM\oplus
V[1] \oplus V^* \oplus \g$, the direct sum of the Eisenberg Lie
algebra  $E=\RM\oplus V[1] \oplus V^*$ and the Lie bialgebra $(\g,[r,-])$. We will now deduce our main result
from:
\begin{proposition}\label{prop:lh}
There exists a $L_\infty$-quasi-isomorphism $\varphi_\hbar$
between 
$(C(\tilde{\g}),[-,-],[r,-])$ and $(T_+\tilde{U},[-,-]_\hbar,[1\otimes 1,-]_\hbar)$. Here $\tilde{U}$ is the Etingof-Kazhdan quantization of $\tilde{\g}$ and so $\tilde{U}=U(E)\otimes U_\hbar(\g)$ where
$U_\hbar(\g)$ is the Etingof-Kazhdan quantization of $(\g,[r,-])$.
The bracket $[-,-]_\hbar$ denotes the Gerstenhaber bracket constructed
in Section \ref{sec:reso} corresponding to the coproduct of $\tilde{U}$.
\end{proposition}

\medskip

\section{Deformed structures and local $L_\infty$-morphism}\label{sec:phi}

\medskip

\subsection{Deformed structures}

Suppose we are
given $\Phi \in (U((\g)^{\otimes 3})^\g[[\hbar]]$ an associator. In particular,
$\gamma^{\otimes 3}(\Phi)$ commutes
with all the invariant differential operator.
We have in fact:
$$[C,\gamma^{\otimes 3}(\Phi)]_G=0\hbox{ for all }C \in U(\tilde{\g})[[\hbar]]).$$
From now on, we will consider the
tensor category of 
$(U(\g)[[\hbar]], \Phi)$-modules (in which 
we want to construct an associative star-product). 
Let us define the ``deformed'' Gerstenhaber bracket
as the Bracket defined in Section \ref{sec:reso} but in the
new tensor category.
We get a new Lie algebra structure on 
$U(\tilde{\g})[[\hbar]]$
given by the bracket $[-,-]_\Phi$ 
defined, for $D,E \in U(\tilde{\g})[[\hbar]]$, by 
$$[D,E]_\Phi=\{D | E\}_\Phi - {(-1)}^{|E||D|}\{E| D\}_\Phi,$$
where for $D \in U(\tilde{\g})^{\otimes d}$ and $E \in U(\tilde{\g})^{\otimes e}$,
$$\{D|E\}=\sum_{i\geq
0}{(-1)}^{(e-1){\cdot}i}\tilde{\Phi}
D^{1,\dots,i,i+1\cdots i+e,i+e+1,\cdots}E^{i+1,\dots,i+e}.$$
$\!\tilde{\Phi}$ corresponds to the obvious change of parenthesis
in the tensor category of 
$(U(\g)[[\hbar]],\!\Phi)$-modules.
For example, if $A$ and $B$ are two $2$-cochains in $U(\tilde{\g})[[\hbar]]$,
one has
$$\{A,B\}_\Phi=A^{12,3}B^{1,2}-\Phi^{-1} A^{1,23}B^{2,3}.$$
\begin{remark}
One could also define a deformed bialgebra structure on 
$(\widehat{\cT} T_+U(\tilde{g})[[\hbar]])^\vee$ and so using Etingof-Kazhdan
dequantization a $G_\infty$-structure on $U(\tilde{\g})[[\hbar]]$ (proof can
be copied from \cite{Ta} or \cite{GH}).
\end{remark}

\medskip

\subsection{Twist quantization of coboundary Lie bialgebras}

\medskip

Let us recall results from \cite{EH}: 
\begin{theorem}\cite{Ha}
Let $(\a,[r,-])$ be a coboundary Lie
bialgebra. There exists a coboundary quantization of it: $(U_\hbar(\a),\Delta_\hbar,R_\hbar)$.
\end{theorem}
\noindent Then, following \cite{Dr}, it was proved in \cite{EH}:
\begin{theorem}
There exists a deformation $\a_\hbar$ of $\a$ in the category of topologically free $\RM[[\hbar]]$-Lie algebras, $J = 1+\hbar r/2+O(\hbar^2)\in U(\a_\hbar)^{\otimes 2}$ and $\Phi_0 \in (U(\a_\hbar)^{\otimes 3})^{\a_\hbar}$ such that the coboundary Hopf algebra $(U_\hbar(\a),\Delta_\hbar,R_\hbar)$ is twist equivalent through $J$ to the coboundary quasi-Hopf algebra $(U(\a_\hbar),\Delta_0,1,\Phi_0)$ and we have, in $U_\hbar(\a)$,
$$J\Delta_0 J^{-1}=\Delta_\hbar\hbox{ and }J^{1,2}J^{12,3}=J^{2,3}J^{1,23}\Phi_0.$$
\end{theorem}
\medskip

\subsection{local $L_\infty$-morphism}

\medskip

Let us keep the notation of the previous section for $\a=\tilde{g}$,
the Lie algebra define in Section \ref{S5.4}. Let us set $F=J^{-1}$ and
$\Phi=\Phi_0$.
We can now prove the existence of a $L_\infty$-morphism for our structures,
in the local case:
\begin{theorem}
There exists a $L_\infty$-quasi-isomorphism $\varphi_\loc$ between the differential Lie algebra
$(\widehat{S}(V)\otimes \Lambda(V^*\oplus \g),[-,-], [r,-])$ (corresponding to local
invariant polyvector fields) and the Lie algebra
$(T_{\widehat{S}(V)}(U(V^*\oplus\g)[[\hbar]]), [-,-]_\Phi,[F,-]_\Phi)$ (corresponding
to invariant polydifferential operators).
\end{theorem}
\dem
Let us consider the Lie bialgebra $\tilde{g}$ defined in Section \ref{sec:Ha}.
Let $\tilde{U}$ be its Etingof-Kazhdan quantization.
We know from Section \ref{sec:Ha} that there exists a $L_\infty$-quasi-isomorphism $\varphi_\hbar$: $(C(\tilde{\g}),[-,-],[r,-])=(\widehat{S}(V)\otimes \Lambda(V^*\oplus \g),[-,-], [r,-])\to (T_+\tilde{U},[-,-]_\hbar,[1\otimes 1,-]_\hbar)$.

Let now define $\varphi_F$: $T_+U_\hbar(\tilde{g})\to T_+U(\tilde{g}_\hbar)[[\hbar]]$ to be the map defined as follows: for $x \in U_\hbar(\tilde{g})^{\otimes n}$, 
$$\varphi_F(x)= F^{12\cdots n-1,n} \cdots F^{12,3}F^{1,2} \cdot x,$$
in $U(\tilde{g}_\hbar)^{\otimes n}[[\hbar]]$. It is clear that $\varphi_F$ is
an isomorphism of differential Lie algebras sending the bracket $[-,-]_\hbar$
to $[-,-]_\Phi$ and the differential $[1\otimes 1,-]_\hbar$ to
$[F,-]_\Phi$. Composing $\varphi_\hbar$ with $\varphi_F$ we get a $L_\infty$-quasi-isomorphism:
$$(\widehat{S}(V)\otimes \Lambda(V^*\oplus \g),[-,-], [r,-])\to (T_+U(\tilde{g}_\hbar),[-,-]_\Phi,[F,-]_\Phi).$$
This gives the result as one can identify $(T_+U(\tilde{g}_\hbar),[-,-]_\Phi,[F,-]_\Phi)$ with $(T_{\widehat{S}(V)}(U(V^*\oplus\g)[[\hbar]]), [-,-]_\Phi,[F,-]_\Phi)$ as differential Lie algebras.
\findem

\smallskip

\begin{remark}
Construction of the Lie algebras isomorphism $\varphi_F$ can be
generalized to differential Lie algebras between any two twist equivalent quasi-Hopf algebras $(H_1,\Delta_1,\Phi_1)$ and $(H_1,\Delta_1,\Phi_1)$
as far as the associators $\Phi_1$ and $\Phi_2$ are invariant (so that one can define corresponding Lie algebras on $T_+H_i$). 
\end{remark}

\medskip

\section{Globalization and Proof of Theorem \ref{theoprinc}}\label{sec:glob}

\medskip

\subsection{Globalization}

\medskip

In this section $X$ is a principal $G$-bundle over a manifold $M$.
We will use Kontsevich globalization procedure as described in \cite{Do}.
One can deduce global version of the local formality theorem (proved in
Section \ref{sec:phi}) to a global one using Fedosov resolution as described
in \cite{Do}. The only things one has to check are the extra conditions
that the $L_\infty$-quasi-isomorphism $\varphi_\loc$ has to fulfill:
\begin{enumerate}\label{eq:glob}
\item The $L_\infty$-quasi-isomorphism $\varphi_\loc$ is equivariant with respect to linear transformations of coordinates.
\item 
$\varphi(v_1,v_2)=0$
for any formal vector fields $v_1$ and $v_2$.
\item If $n \geq 2$ and
$v$ is a linear vector filed in the coordinates on 
$\RM^n$, then for any set of polyvector fields
$\gamma_2,\dots,\gamma_n$ we have
$\varphi(v,\gamma_2,\dots, \ve \gamma_n)=0.$
\end{enumerate}
\begin{proposition}\label{pro-cnd}
The $L_\infty$-quasi-isomorphism $\varphi_\loc$ can be built so that it satisfies 
those three conditions
\end{proposition}
\dem
Let us recall that the map $\varphi_\loc$ was built from two differential Lie algebra morphism:
$\varphi_{\Ger_\infty}$: $C(\tilde{g}_\hbar) \to \CuTTUt$ and
$\varphi_{\G_\infty}$: $T_+\tilde{U} \to \CuTTUt$. Recall also that $\varphi_{\G_\infty}$ was
built as a resolution of $T_+\tilde{U}$.
Now, instead of using Theorem \ref{theo:koninv} we will construct directly the
composition of $\varphi_{\Ger_\infty}$ with the inverse of $\varphi_{\G_\infty}$.
More precisely, we will construct $\tilde{\varphi}_{\Ger_\infty}$ a $L_\infty$-quasi-isomorphism deforming
$\varphi_{\Ger_\infty}$ so that the image of $\tilde{\varphi}_{\Ger_\infty}$ is
contained in the space of cocycle of $\CuTTUt$, and satisfy conditions of the theorem. This will then give the result.
Let us now recall a useful lemma that can be found in \cite{Do}:
\begin{lemma}
Let $\phi$ be a $L_\infty$-quasi-isomorphism between two DGLAs $(\g_1,d_1)$ and $(\g_2,d_2)$.
Let $\phi^{n}$: $S^n(\g_1[1])\to \g_2$ be the structure
maps of $\phi$. Let $m\geq 1$. Then it is possible to construct a deformed
$L_\infty$-quasi-isomorphism $\tilde{\phi}$ satisfying
\begin{itemize}
\item $\tilde{\phi}(\gamma_1,\dots,\gamma_n)=\phi(\gamma_1,\dots,\gamma_n)$, for
$n<m$.
\item $\tilde{\phi}(\gamma_1,\dots,\gamma_m)=\phi(\gamma_1,\dots,\gamma_m)
+ d_2V(\gamma_1,\dots,\gamma_m)\\
\hbox{}\hskip 2cm -\sum_{1\leq l \leq m}{(-1)}^{l+|\gamma_1|+\cdots+|\gamma_{l-1}|}V(\gamma_1,\dots,d_1\gamma_l,\dots,\gamma_m)$,\\
where $V$: $S^m(\g_1[1])\to \g_2$ is an arbitrary polylinear map.
\end{itemize}
\end{lemma}
Moreover one has explicit computation of $\tilde{\phi}$ from $\phi$ and $V$:
let $D_1=d_1+d_1^{1,1}$ and $D_2=d_2+d_2^{1,1}$ be the structure maps of $\g_1$ and $\g_2$ and $\Delta_1, \Delta_2$ the associated free comultiplications (see Section \ref{sec:structure}). Then, for $x\in C(\g_1[1])$,
$$\tilde{\phi}(x)=\phi(x)+D_2V(x) +V(D_1x),$$
where $V$ is extended as follows:
$$\Delta_2D_2V(x)=\left(\phi \otimes V + V \otimes \phi+\frac{1}{2}(V \otimes D_2V+D_2V\otimes V)+\frac{1}{2}(V \otimes VD_1+VD_1\otimes V)\right)\Delta_1(x).$$
Let us denote by $\partial$ the differential in $\CuTTUt$
and $\phi$ will be the map $\varphi_{\Ger_\infty}$. We know that the complex $(\CuTTUt,\partial)$ is acyclic except for elements of $T_+\tilde{U}$. We will write
$(x)_0$ for the component in $T_+\tilde{U}$ of  an element $x\in \CuTTUt$.
So, as $\partial \phi^1(x)=0$, there exists a linear map $V$: $C(\tilde{g}_\hbar) \to \CuTTUt$ such that for every element $v \in C(\tilde{g}_\hbar)$, $\phi^1(x)=(\phi^1(x))_0+\partial V(x)$. Note that for degree reasons, when $x$ is a vector field, $(V(x))_0$ is
a function and so can be chosen to be zero. Moreover, the map
$\phi$ is equivariant with respect to change of coordinates (see \cite{Ha2}, Section 5.2). So we can assume that $V$ is also equivariant.
So one can define $\tilde{\phi}$ as in the lemma and $\tilde{\phi}$
satisfies the first condition of the theorem. Moreover $(\tilde{\phi})_0$ clearly satisfy the second condition and the third is again a consequence of equivariance with respect to change of coordinates. Let us replace $\phi$ with $\tilde{\phi}$.
We will know proceed by induction and suppose that the first condition of the theorem is true the structure maps of $\phi$, that the second and third conditions a true for their $(-)_0$ parts,
and are also true for the structure maps of $\phi^i=(\phi^i)_0$ for $i\leq n$. Using the induction hypothesis and the fact that $\phi$ is a $L_\infty$-morphism, we get that $(\partial \phi^{n+1})_0=0$.
So there exists $W$ such that $\phi^{n+1}=(\phi^{n+1})_0 + \partial W$. Again $W$ can be chosen equivariant. $(W)_0$ can be chosen to satisfy the last condition (the second is automatic for degree reason) as $(\phi^{n+1})_0$
satisfies it. Then again thanks to equivariance, one checks that $\tilde{\phi}$
obtained from $\phi$ and $W$ satisfies the hypothesis of the induction. 
This concludes the proof.
\findem

\medskip

\subsection{Proof of Theorem \ref{theoprinc}}

\medskip

Let us summarize what we have done so far:
\begin{theorem}\label{theo:linfprinc}
Let $X$ be a principal $G$-bundle over a manifold $M$. Let $r \wedge^2 \g$ such that 
$[r,r]=Z \in (\wedge^3 \g)^\g$. There exists $\Phi=1+\frac{\hbar^2}{6} Z + O(\hbar^3) 
\in (U(\g)^{\otimes 3})^\g$, $J=1+\hbar r + O(\hbar^2) 
\in U(\g)^{\otimes 2}
[[\hbar]]$ such that $J^{1,2}J^{12,3}=J^{2,3}J^{1,23}\Phi$, a deformation $\g_\hbar$ of the
Lie algebra $\g$ and $\varphi$ a $L_\infty$-quasi-isomorphism
$$\varphi:~(\Tpi(M),[-,-]_S,[r,-]_S) \to (\Dpi(M),[-,-]_\Phi,[J,-]_\Phi),$$
where $\Tpi(M)$ and $\Dpi(M)$ are respectively the spaces of invariant polyvectors fields an $M$ or polydifferential operators on $M$ and  $[-,-]_\phi$ is the deformed Gerstenhaber bracket in the tensor category of
$(U(\g_\hbar)[[\hbar]],\Phi)$-modules.
\end{theorem}

Suppose now that $(X,\pi,r,Z)$ is a quasi-$(r,Z)$-Poisson manifold.
Set $\pi=\pi'+r$. Then $\pi$ is a Maurer-Cartan in the DGLA
$(\Tpi(M),[-,-]_S,[r,-]_S)$. Thanks to Theorem \ref{theo:linfprinc}, we know
that those Maurer-Cartan elements, up to Gauge transform are
in one to one correspondence with Maurer-Cartan elements of
the DGLA $(\Dpi(M),[-,-]_\Phi,[J,-]_\Phi)$, up to Gauge transform.
If $m_\star'$ is such a Maurer-Cartan element, set $m_\star=m_\star' + J$,
$m_\star$ is a quantization of the quasi-$(r,Z)$-Poisson manifold $(X,\pi,r,Z)$.
We have prove:
\begin{theorem}
Let $(X,\pi,r,Z)$ be a quasi-$(r,Z)$-Poisson manifold, quantization of $X$,
up to equivalence, are in one to one correspondence with 
$$\{\pi_\hbar= \hbar \pi + O(\hbar^2)\hbox{ such that }[r,\pi_\hbar]+\frac{1}{2}[\pi_\hbar,\pi_\hbar]=0\}.$$
\end{theorem}

\medskip

\section{Quantization of modified dynamical Yang-Baxter $r$-matrices}\label{sec:discus}

\medskip

We know that modified 
classical dynamical Yang-Baxter $r$-matrices provide examples of quasi-Poisson manifolds.
Let us recall their definition:
let $\h$ be a Lie subalgebra of $\g$. Let
$\rho$ be a $\h$-equivariant map $\rho$~:
$\h^* \to \ve^2(\g)$, solution of the modified 
classical dynamical Yang-Baxter equation:
$$-\Alt(\d \rho) + \CYB (\rho)=Z,$$
where
$$\CYB(\rho)=[\rho^{1,2},\rho^{1,3}]+
[\rho^{1,2},\rho^{2,3}]+[\rho^{1,3},\rho^{2,3}]$$
and 
$$\Alt(\d \rho)=\sum_i h_i^1 \frac{\partial \rho^{2,3}}{\partial \lambda_i}-
\sum_i h_i^2 \frac{\partial \rho^{1,3}}{\partial \lambda_i}+ \sum_i h_i^3 \frac{\partial \rho^{1,2}}{\partial \lambda_i}.$$
Using a quasi-Poisson generalization of a construction of Xu \cite{Xu},
Enriquez and Etingof 
\cite{EE1} built a quasi-Poisson manifold $X_\rho$ associated to $\rho$
(for which the action of the corresponding group $G$ is free).
They then prove (following \cite{Xu}) that any twist quantization 
$J$ associated to an associator $\Phi$ 
(i.e. 
\begin{itemize}
\item $J \in \Mer(\h^*,U(\g)^{\otimes2}[[\hbar]]$, $\h$-invariant, such that
$J(\lambda)=1+O(\hbar)$, 
 \item $J^{12,3}(\lambda) \star J^{1,2}(\lambda + \lambda \hbar^3)=
\Phi^{-1} J^{1,23}(\lambda) \star J^{2,3}(\lambda),$
\item $Z=\Alt\left(\frac{\Phi-1}{\hbar^2}\right) \mod \hbar$,
\item $\rho(\lambda)=\left(\frac{J(\lambda)-1}{\hbar}\right)
-\left(\frac{J(\lambda)-1}{\hbar}\right)^{2,1} \mod \hbar$) 
\end{itemize}
gives rise to a quantization
of the quasi-Poisson manifold $X_\rho$. 
Our result provides us with a quantization of the manifolds $X_\rho$ 
when $Z$ satisfies our conditions but
unfortunately, we don't
know whether this quantization provides us with a twist quantization
of the modified dynamical $r$-matrix $\rho$. 

\medskip

Let us write, according to \cite{Xu}, the Poisson bracket associated $\pi_\rho$ to 
a dynamical $r$-matrix $\rho$: in the decomposition of $\Tpi=\wedge T\h^*\otimes
\wedge \g$,
$$\pi_\rho=\pi_{\h^*}+\sum_i\frac{\partial}{\partial \lambda_i}\wedge \h_i + \rho \in \wedge^2T\h^* \oplus T\h^*\wedge \g \oplus \wedge^2\g,$$
where $h_i$ and $\lambda_i$ are basis and dual basis of $\h$ and $\pi_{\h^*}$ is the Kostant-Kirilov-Souriau Poisson bracket
(which we will denote $\star_{\h^*}$). Now, still following \cite{Xu}, 
a quantization $\star$ of $X_\rho$ corresponds to a quantization of $\rho$ is an only if:
it satisfies the following conditions:
\begin{enumerate}
\item for any $f,g \in C^\infty(\h^*),$ $f(\lambda)\star g(\lambda)=f\star_{\h^*}g$,
\item for any $f \in C^\infty(G)$ and $g\in C^\infty(\h^*),$
$f(x)\star g(\lambda)=f(x)g(\lambda)$,
\item for any $f \in C^\infty(\h^*)$ and $g\in C^\infty(G),$
$f(\lambda)\star g(x)=\sum_{k \geq 0}\frac{\hbar^k}{k!}\frac{\partial^k f}{\partial\lambda^{i_1}\cdots \partial\lambda^{i_k}}h_{i_1}\cdots h_{i_k} g$,
\item for any $f, g \in C^\infty(G)$,
$f(x)\star g(x)=R(f,g),$ where $R$ would be the quantization of $\rho$.
\end{enumerate}
Let us notice that our quantization of $X_\rho$ will not satisfiy
those conditions as, for symmetry conditions, conditions 2 and 3 will
not be fulfilled.
So we will use a trick proved by Alekseev and Calaque (\cite{AC}).
Let us first recall their definition of strongly $\g$-invariant quantization of quasi-Poisson manifold.
\begin{definition}
Let $\star$ be a quantization of a quasi-Poisson manifold $(X,\pi,Z)$.
Suppose $\mu$: $X \to \h^*$ is a momentum map for which the map $M=U(\mu^*)\circ$ sym: $(\OC_{\h^*}[[\hbar]],\star_{PBW})\to (\OC_X[[\hbar]],\star)$
is an algebra morphsism satisfying $[M(x),f]_\star=\hbar\{\mu^*x,f\}$
for any $f\in \OC_X$ and any $x \in \g$, then we say that the quantization $\star$ is strongly $\g$-invariant.
\end{definition}
\begin{proposition}\cite{AC}\label{propAC}
Assume that $\star$ is a strongly $\g$-invariant quantization of $X_\rho$. Then there exists a gauge equivalent quantization $\star'$ of $X_\rho$ such that corresponds to a quantization of $\rho$.
\end{proposition}
Thus we only need to prove that we can construct a strongly $\g$-invariant quantization of $X_\rho$.
\begin{theorem}
Suppose that $\h$ is an abelian subalgebra of $\g$ and that there exists a decomposition $\g = \h \oplus \m$ with $[\h,\m]\subset \m$.
Then the modified classical dynamical Yang-Baxter $r$-matrice ($r-Z)$ can be quantized.
\end{theorem}
\dem
Using Proposition \ref{propAC}, we need to prove that there exists a 
a strongly $\g$-invariant quantization of the quasi-Poisson 
manifold $X_\rho$. To do so, we can copy the proof of Proposition
\ref{pro-cnd} to get adapted product.
\findem


\vskip20pt

\smallskip

\begin{thebibliography}{99}

\bibitem[AC]{AC}  A. Alekseev, D. Calaque, 
{\it Quantization of symplectic dynamical r-matrices and the quantum composition formula},
Communications in Mathematical Physics 273 (2007), no. 1, 119--136.

\bibitem[AK]{AK}  A. Alekseev, Y. Kosmann-Schwarzbach, 
{\it Manin pairs and moment maps},
J. Differential Geom. 56 (2000), no. 1, 133--165.

\bibitem[AKM]{AKM}  A. Alekseev, Y. Kosmann-Schwarzbach, E. Meinrenken,
{\it Quasi-Poisson manifolds},
Canad. J. Math. 54 (2002),  no. 1, 3--29.



\bibitem[BFFLS1]{BFFLS1}  F. Bayen, M. Flato, C. Fronsdal, A. ichnerowicz, D.
Sternheimer,
{\it Quantum mechanics as a deformation of classical mechanics},
    Lett. Math. Phys. 1 (1975/77),  no. 6, 521--530.


\bibitem[BFFLS2]{BFFLS2} F. Bayen, M. Flato, C. Fronsdal, A. Lichnerowicz, D.
Sternheimer, 
{\it  Deformation theory and quantization, I and II}, 
Ann. Phys. 111 (1977), 61-151.

\bibitem[Do]{Do} V. Dolgushev,
{\it Covariant and equivariant formality theorems}, Adv. Math.  {\bf 191}  (2005),  no. 1, 147--177.

\bibitem[Dr]{Dr}
V. Drinfeld,
``Quasi-Hopf algebras". Leningrad Math. J. {\bf 1} (1990),
no. 6, 1419--1457.

\bibitem[EE1]{EE1} B. Enriquez, P. Etingof,
{\it  Quantization of Alekseev-Meinrenken dynamical r-matrices},
math.QA/0302067.

\bibitem[EH]{EH} B. Enriquez, G. Halbout, {\it Quantization of coboundary Lie bialgebras}, to appear in Annals of Math.

  
\bibitem[EK]{EK} P. Etingof, D. Kazhdan,
  {\it  Quantization of Lie bialgebras. I}, Selecta Math. 
  (N.S.) 2 (1996), no.
  1, 1--41
    

    


\bibitem[Fr]{Fr} C. Fronsdal, {\it Some ideas about quantization},
Rep. Math. Phys.  15  (1979), no. 1, 111--145.

\bibitem[Gav]{Gav} F. Gavarini, {\it The quantum duality principle,}
Ann. Inst. Fourier (Grenoble), {\bf 52} (2002), no. 3, 809-834.

  
  

  

\bibitem[GH]{GH} G. Ginot, G. Halbout,
{\it A formality theorem for Poisson manifolds},
 Lett. Math. Phys.  66  (2003),  no. 1-2.


  
\bibitem[Ha]{Ha} G. Halbout, {\it Formality theorem for Lie bialgebras and quantization of twists and coboundary $r$-matrices}.  Adv. Math.  207  (2006),  no. 2, 617--633.

\bibitem[Ha2]{Ha2} G. Halbout, {\it from associators to a global formulation.  Ann. Math. Blaise Pascal}  13  (2006),  no. 2, 313--348.



\bibitem[Ko]{Ko} M. Kontsevich, {\em Deformation quantization of Poisson manifolds}.  Lett. Math. Phys.  {\bf 66}  (2003),  no. 3, 157--216.
 

\bibitem[Me]{Me} S. Merkulov, {it Quantization of strongly homotopy Lie bialgebras}, math.QA/0612431

\bibitem[Ta]{Ta} D. Tamarkin,
  {\it Another proof of M. Kontsevich's formality theorem}, math.QA/9803025.

\bibitem[Xu]{Xu} P. Xu, {\it Quantum dynamical Yang-Baxter equation
over a nonabelian base}, Commun. Math. Phys. 226 (2002), no. 3, 475--495.



\end{thebibliography}
\enddocument